\documentstyle[12pt]{article}
 
\def\a{\alpha}\def\b{\beta}\def\o{\omega}\def\l{\lambda}
\def\U{\Upsilon}\def\P{\hbox{Pr}}\def\ch{\choose}
\begin{document}
 
\thispagestyle{empty}
\title{\bf Efficient Pooling Designs for Library Screening}
\nopagebreak \author{ W. J. Bruno$^{1,2}$, E. Knill$^{1,3}$, D. J. Balding$^5$,\\
                D. C. Bruce$^{1,4}$, N. A. Doggett$^{1,4}$, W. W. Sawhill$^{1,4}$,\\
                R. L. Stallings$^6$, C. C. Whittaker$^{1,2}$, and D. C. Torney$^{1,2,7}$\\}
\maketitle
 
\begin{center}
$^1$ Center for Human Genome Studies\\ Los Alamos National Laboratory
\\ Los Alamos, NM 87545 \medskip
 
$^2$ T-10, Theoretical Biology and Biophysics, Mailstop K710 \\
Theoretical Division \medskip
 
$^3$ CIC-3, Computer Research and Applications, Mailstop K990 \\
Computing, Information, and Communications Division \medskip
 
$^4$ LS-2, Genomics and Structural Biology, Mailstop M880\\
Life Sciences Division \medskip
 
$^5$ School of Mathematical Sciences \\ Queen Mary and Westfield
College \\ University of London \\ London, E1 4NS, UK \medskip
 
$^6$ Department of Human Genetics \\ University of Pittsburgh \\
Pittsburgh, PA, 15261, USA \medskip
 
$^7$ Author for Correspondence \\
     dct@ipmati.lanl.gov
\end{center}
\newpage
 
\begin{abstract}
 
We describe efficient methods for screening clone libraries, based on
pooling schemes which we call ``random $k$-sets designs''.  In these
designs, the pools in which any clone occurs are equally likely to be
any possible selection of $k$ from the $v$ pools. The values of $k$
and $v$ can be chosen to optimize desirable properties.  Random
$k$-sets designs have substantial advantages over alternative pooling
schemes: they are efficient, flexible, easy to specify, require fewer
pools, and have error-correcting and error-detecting capabilities.  In
addition, screening can often be achieved in only one pass, thus
facilitating automation.  For design comparison, we assume a binomial
distribution for the number of ``positive'' clones, with parameters
$n$, the number of clones, and $c$, the coverage.  We propose the
expected number of {\em resolved positive} clones---clones which are
definitely positive based upon the pool assays---as a criterion for
the efficiency of a pooling design.  We determine the value of $k$
which is optimal, with respect to this criterion, as a function of
$v$, $n$ and $c$.  We also describe superior $k$-sets designs called
$k$-sets packing designs.  As an illustration, we discuss a
robotically implemented design for a 2.5-fold-coverage, human
chromosome 16 YAC library of $n=1,298$ clones.  We also estimate the
probability each clone is positive, given the pool-assay data and a
model for experimental errors.
 
\end{abstract}
\setlength{\baselineskip} {0.35in}
\section{Introduction}
 
Much of the current effort of the Human Genome Project involves the
screening of large recombinant DNA libraries in order to isolate
clones containing a particular DNA sequence.  This screening is
important for disease-gene mapping and also for large-scale clone
mapping \cite{ols}. More generally, efficient screening techniques can
facilitate a broad range of basic and applied biological research.
Whenever the objective is to find ``needles in a haystack'', a
reliable test indicating whether or not at least one needle occurs in
a specific part of the haystack can greatly facilitate the isolation
of the needles \cite{du}, \cite{dy}.  Such tests are called ``binary
group tests''.  The most reliable binary group tests routinely used to
screen groups of clones for a particular DNA sequence employ either a
hybridization-based assay or a PCR-based, STS assay.
 
Each group of clones is called a ``pool'', a group test is called a
``pool assay'', and a collection of pools is called a ``pooling
design''.  A convenient specification of a pooling design is a binary
clone-{\em by}-pool incidence matrix: if a clone occurs in a pool the
matrix element equals unity.  Figure 1a depicts an incidence matrix
for a small pooling design.
 
Any clone containing a DNA sequence of interest is termed a
``positive''.  Pools yielding a positive assay are termed ``positive''
pools.  For the time being we assume there are no experimental errors,
in which case each positive pool contains one or more positive clones.
After the pools have been assayed, the clones fall into one of four
conceptual categories, illustrated in Figure 1b. If a negative clone
occurs in at least one pool containing only negative clones then its
negative status can be determined from the pool assays and it is
called a ``resolved negative''.  The remaining negative clones are
called an ``unresolved negative''.  A positive clone is called a ``resolved
positive'' if it occurs in at least one pool with no other positive
clone and no unresolved negative clone, otherwise it is called
``unresolved positive''.  Thus, the status of the resolved positive
and resolved negative clones is resolved by the experiment.  Although
the experimenter will in principle be unable to distinguish the
unresolved positive from the unresolved negative clones, as will be
seen, it may be useful to separately analyze these two types of
clones.
 
The expected number of unresolved negative clones was previously
proposed as a criterion for selecting pooling designs \cite{bar}.
However, the library-screening objective can range from the isolation
of a small proportion of the positive clones to the isolation of all
of the positive clones---hence, we propose the expected number of
resolved positive clones as a criterion of design optimality.  For
example, a cDNA library typically contains multiple copies of many
sequences, but the experimenter may want to identify only one of these clones.
In Section \ref{res} we calculate both these expectations for several
proposed and implemented pooling designs, under the assumption of a
binomially-distributed number of positive clones.  Confirmatory tests
are ordinarily performed on the candidate positive clones.  The
average number of such tests that would be required to determine the
status of all the clones equals the expected number of positive clones
plus the expected number of unresolved negative clones.  The
probability that a design yields a one-pass solution (i.e.\ the status
of all clones is resolved) may also provide a useful criterion for
design optimality.  This criterion, introduced in \cite{b&t}, will be
discussed in a future publication \cite{ima}.
 
Our focus is on one-stage designs in which all of the pools are
assayed in one pass. The use of one-stage designs facilitates
automation because a robot can be fully programmed from the outset,
irrespective of intermediate results. In addition, one-stage designs
require the construction of fewer pools over multiple screenings,
since the same pools are used for each screening.
Standard pooling designs for large libraries have
several stages, and they often involve row-and-column pools, {\em
viz.}\ Figure 2 and \cite{ame}, \cite{chu}, \cite{eva}, \cite{gre},
\cite{slo}.  The examples of Section \ref{comp} illustrate the
differences between these approaches for a hypothetical
tenfold-coverage, human-genome library with 72,000 clones.

We propose ``random $k$-sets'' designs for one-stage library
screening.  In these designs, each clone occurs in $k$ pools, and all
choices of the $k$ pools are equally likely \cite{dy1}.  Random
$k$-sets designs are easy to specify for any number of pools, and they
are efficient, in terms of the expected number of resolved positives,
in comparison to alternative designs.  The expected numbers of
resolved negative and resolved positive clones are given in Section
\ref{mm}.  Using these formulas, the optimal choices of $k$ for a
given application can be determined.  Numerical results are given for
a range of clone library sizes and coverages in Section \ref{res}.
(Although performance will vary according to the specific instance of
the random design, measures of performance are narrowly distributed
between different instances for large designs; thus, expected measures
of performance provide a useful guide).  We also discuss techniques
for constructing superior $k$-sets designs, called ``$k$-sets
packing'' designs, in Section
\ref{gpd}.
 
In Section \ref{ee} we describe techniques for ranking candidate
positive clones, given the pooling design, the pool-assay data, and a
model for experimental errors.  Ranking methods are employed to
illustrate the performance of random $k$-sets designs in the presence
of pool-assay errors, through computer simulation.
 
We have previously described efficient techniques for manual
construction of pools \cite{mkm1}.  The designs proposed here are
intended for implementation using robots.  We employed the Packard 204
multiPROBE robot to pool a 1,298-clone, 2.5-fold-coverage,
chromosome-specific YAC library, implementing a four-sets design on 47
pools.
 
\subsection{Notation}
 
The number of clones in the library is denoted by $n$.  The number of
pools in a pooling design is denoted by $v$, and the number of pools
in which a clone occurs is denoted by $k$.  The number of clones in a
clone library which cover any particular location of the cloned region
is assumed to be binomially-distributed with expectation $c$ 
(called the coverage).
The standard notation Pr$(A|Z)$ is used for the conditional
probability of event $A$, given event $Z$.  We use the notation ${a\ch
b}$ for the binomial coefficient $a!/\left(b!(a-b)!\right)$, and we
introduce $B$ for the binomial probabilities
\begin{equation}
B(a,b,t)={a\ch b}t^b(1-t)^{a-b}.
\label{bin}
\end{equation}
 
\section{Methods}
\label{meth}
 
We developed several computer programs as part of our pooling
methodology, described below.  Copies of these programs are freely
available from the author for correspondence.
 
\subsection{Pooling and Screening Methods for a 1,298-Clone YAC Library}
\label{imp}
 
A complete-digest, Cla I library of 1,298 clones was generated from
human chromosome 16 \cite{mkm}.  Based on an average cloned-DNA size
of 200 kilobasepairs and a chromosome size of 98 megabasepairs
\cite{mor}, the coverage of this library was conservatively estimated
at $c=2.5$.  A four-sets packing design on $v=47$ pools was used; see
Section \ref{gpd} for its characterization.  Thus, each clone occurred
in four of the pools.  (With this value of $v$, two sets of pools and
controls can be put on one 96-well dish).  Deep-96-well mitrotiter
dishes were filled with YPD medium using the Dynatech plate-filler;
each well contained 1.9 ml of medium.  These dishes were incubated at
$30^{\circ}$ C for 48 hours to check for contamination prior to
inoculation.  Then, 14 uncontaminated dishes were inoculated by hand
stamping from a selective-medium copy of the library, lacking uracil
and tryptophan.  The inoculated dishes were incubated at $30^{\circ}$
C for 72 hours, achieving stationary growth.  We subsequently used a
Packard multiPROBE 204 robot to construct the pools.  The robot cycled
through the following procedure, using all four liquid-handling tips:
suspend the yeast in a well, aspirate 0.8 ml, dispense 0.4 ml in each
of two specified pool tubes, and sterilize and clean the tips.  (A
larger syringe would have allowed us to aspirate the full well volume
and would have reduced the overall time for pooling).  The sterilizing
and tip cleaning began with a five-second rinsing with five ml of
distilled water.  Next, one ml of 0.525\% sodium hypochlorite solution
was aspirated and discharged.  Then, the tips were rinsed with five ml
of sterile water, and the tips were placed in individual cups to clean
the outside as well.  (These rinsing and sterilizing steps account for
about half of the run time.  They were sufficient to prevent
cross-contamination on the basis of a yeast-growth assay).  Final pool
volumes of approximately 46 ml were collected in 50 ml centrifuge
tubes.  Pools were maintained at room temperature during the pooling
run, which took approximately 12 hours.  Human intervention was
required hourly to replace the dishes from the arrayed library.
 
Agarose gel plugs containing yeast cells were prepared from each pool.
The plugs were treated with zymolase to make spheroplasts, and they
were treated with ESP (500 mM EDTS, 1\% sarcosyl, and 1
mg/ml proteinase K), to digest away proteins.  The plugs were dialyzed
extensively against TE.  The DNA was purified from the agarose using
the Geneclean kit.  GeneAmp PCR reagent kits were used for PCR
screening of the pool DNA.  (50 ng of pool DNA and 2.5 units of Taq
polymerase and 1.5 mM MgCl$_{2}$ were the noteworthy features of the
50 $\mu$l reaction volume).
 
\subsection{Methods for Programming Robots}
\label{robie}
 
Two software systems, written in `C', were developed to facilitate use
of the Packard multiPROBE 204 robot.  These systems were both used in
the implementation described in Section \ref{imp}.
 
Our ``scheduling'' software can readily be adapted for an arbitrary
robot.  This system enables the implementation of arbitrary pooling
designs, {\em viz.}\ Figure 1a, creating a list of volumes for
well-to-pool transfers.  In addition, it lists the location, size, and
quantity of ``source'' and ``destination'' plastic-ware.  The
scheduling system can be used to coordinate the simultaneous use of
multiple robots engaged in constructing a set of pools from one
library.  The output of the scheduling system is used by a robot
instruction ``interpreter'', described next.
 
The interpreter program reads the scheduling input file and generates
a complete list of commands for the Packard robot.  The basic actions
of the robot are aspirating, dispensing, and rinsing.  Each of these
actions requires many steps and hence many robot commands.  When the scheduling
program
provides the necessary parameters, the interpreter program generates
these commands, which can then be implemented using HAM, the Packard-robot
software controlling the monitor and the robot.
 
For example, the scheduling program supplies the well number and the
volume of aspiration, and the interpreter program generates the
appropriate list of robot commands.  The interpreter determines
whether multiple tips can aspirate or dispense simultaneously.  It
also records all commands and actions and enables termination of the
run at any time.
 
\subsection{$k$-Sets Design Generation---$k$-Sets Packings}
\label{gpd}
 
A four-sets design on 47 pools was constructed using a pseudo-random
number generator.  The design was based on a random four-sets design;
however, a prospective four-set was discarded if it had more than two
elements in common with another four-set already in the design.  In
addition, we ensured that each pool contained approximately the same
number of clones (109 to 111).  In Section \ref{eufour} our
performance criterion is used to compare the design we generated using these
constraints to the average perfomance of a random four-sets design.
 
Designs with a bound on the size of the intersection between two
$k$-sets are referred to as ``$k$-sets packing'' designs. In general,
it is unsatisfactory to have two identical $k$-sets in a random
$k$-sets design: if one of them is positive it would be impossible to
establish the status of the other.  Similarly, it is undesirable for
two clones to coincide in a large number of pools (i.e.\ for the sets
of pools each occurs in to be similar).  Bounding the intersection
ensures that a negative clone cannot be unresolved due to a small
number of positive clones.

\subsection{Mathematical Methods}
\label{mm}
 
In this section we give expressions for the measures of design
performance---expected numbers of resolved negative and resolved
positive clones---for random $k$-sets designs. We assume that each
clone is positive with probability $c/n$, independently of the other
clones, which yields a binomial distribution for the number of
positive clones. Although this may not be strictly realistic, it
provides a reasonable approximation for the purposes of design
comparison \cite{c&c}, \cite{l&w}.  The expressions we derive are
readily modified to accommodate an arbitrary probability that there
are $p$ positive clones, provided that all $p$-subsets of the clones
are equally likely to be the positive clones.  Such a modification
would be useful when screening cDNA libraries, for example.  Here we assume no
experimental error.  Pool-assay errors are treated in Section
\ref{ee}.
 
Recall the definition of random $k$-sets designs: each clone occurs in
precisely $k$ of the $v$ pools, and each of the possible ${v\ch k}$
subsets of $k$ pools has equal probability to be the pools in which a
particular clone will occur.
 
We begin by formulating $\bar{N}$, the expected number of
unresolved negative clones---negative clones occurring only in
positive pools; thus, $n{-}c{-}\bar{N}$ is the expected number of
resolved negative clones.
 
\begin{equation}
\bar{N} = \sum_{p=0}^{n}(n{-}p)B(n,p,c/n)\sum_{i=0}^{k} {k\ch i}(-1)^iz_i^p,
\label{nuexact}
\end{equation}
where $$z_i = {v{-}i\ch k}\,{v\ch k}^{-1},$$ and $B$ denotes the
binomial probabilities (\ref{bin}).  Interchanging the order of
summation in (\ref{nuexact}) and discarding some terms of order
$n^{-1}$ yields
\begin{equation}
\bar{N}\approx(n{-}c)\sum_{i=0}^{k} {k\ch i}(-1)^{i}{\rm e}^{-c(1-z_{i})}.
\label{nuasymp}
\end{equation}
Insight into the behavior of (\ref{nuexact}) can be obtained from the
``independent pools'' approximation. The inner summation is an
inclusion-exclusion formula for the probability, $K^{(p)}_k$, that all
of the $k$ pools in which a negative clone occurs contain at least one
positive clone, given that there are $p$ positive clones.  The
probability that a given pool is negative is $(1{-}k/v)^p$, and the
independent pools approximation uses
\begin{equation}
K^{(p)}_k\approx\left(1-\left(1{-}\frac{k}{v}\right)^p\right)^k.
\label{appr}
\end{equation}
The approximation (\ref{appr}) would be exact if pool outcomes were
independent.  Because pool outcomes are negatively
correlated, the independent pools approximation gives an upper bound
for $K^{(p)}_k$ and hence $\bar{N}$.  The bound is tight for $p$ not
too small and improves as $v$ increases. A recursive expression for
$K^{(p)}_k$ is given in Appendix
\ref{app} (Equation (\ref{covjr})).
 
We now formulate $\bar{P}$, the expected number of unresolved positive
clones; thus, $c{-}\bar{P}$ is the expected number of resolved
positive clones. The details of the calculation are given in Appendix
\ref{app}.  Recall that, in the absence of pool-assay errors, a
positive clone is unresolved if it occurs only in pools containing
either another positive clone or an unresolved negative clone.
 
After substitution and performing the summations over $p$ and $u$,
Equation (\ref{last}) becomes
\begin{equation}
\bar{P} = c\left(1-\frac{c}{n}\right)^{n-1}\a^{(1)}+
c\sum_{x=k}^v\sum_{y=x}^v{v\ch x}
{x\ch y{-}k}{v{-}x \ch y{-}x}{v \ch k}^{-1}Q_{x,y}
\label{rhof}
\end{equation}
in which
$$Q_{x,y}=\sum_{i=v-x}^{v-k}{x\ch v{-}i}(-1)^{i-v+x}
\sum_{j=0}^{y-x}{y{-}x\ch j} (-1)^{j}R_j,$$
$$R_j=\left(\left(\frac{z_ic}{n}\right)+\left(1{-}\frac{c}{n}\right)
(1-\b(1{-}\xi_j))\right)^{n-1}-
\left(\left(1{-}\frac{c}{n}\right)(1{-}\b(1{-}\xi_{j}))\right)^{n-1},$$
and $\a^{(1)}$, $\b$ and $\xi_j$ are defined at, respectively,
(\ref{p1}), (\ref{s}) and (\ref{prguxyp}), and $z_{i}$ is defined
below Equation (\ref{nuexact}).
 
It is seen that (\ref{rhof}) can be evaluated with ${\cal O}(v^{4})$
summands.  To evaluate (\ref{nuexact}) and (\ref{rhof}) for large $v$,
we used MAPLE V's extended-floating-point precision capabilities
\cite{map}.  In addition, we wrote our own extended-precision FORTRAN
routines, which typically required about one percent as much c.p.u.\
time to obtain an equally accurate result.  The evaluation can also be
done without using extended precision: the alternating signs can be
eliminated by using the recursive Equations (\ref{covjr}) and
(\ref{uncovjr}).  In this case the summations over $p$ and $u$ must be
approximated numerically.
 
We have also derived approximations for $\bar{P}$ that are more easily
evaluated.  If the number of positive pools were fixed at the expected
number, $\o$:
 
       $$\o = v\left(1-\left(1{-}\frac{k}{v}\right)^{p}\right),$$
 
then the probability a negative clone is unresolved negative is $\mu$,
given by
 
                $$\mu={\o\ch k}{v\ch k}^{-1},$$
 
Averaging over the number of unresolved negatives, with a Poisson approximation
to the binomial distribution, and using the ``independent pools'' approximation 
gives
 
         $$\bar{P}\approx\sum_{p=1}^{n}pB(n,p,c/n)
\left(1-\left(1{-}\frac{k}{v}\right)^{p-1}{\rm e}^{{-}\mu (n{-}p) k/\o}\right)^{k},$$
 
\noindent
A more accurate approximate formula is obtained if pools are not
assumed to be independent:
 
\begin{equation}
   \bar{P}\approx \sum_{p=1}^{n} pB(n,p,c/n)
                  \sum_{i=0}^{k}{k \ch i}(-1)^{i} z_{i}^{p{-}1}{\rm e}^{{-}\mu(n{-}p)(1-\zeta_i)},
\label{rhoruf}
\end{equation}
 
\noindent
in which
 
          $$\zeta_{i} = {\o{-}i\ch k}\,{\o\ch k}^{{-}1}.$$
 
\noindent
Equation (\ref{rhoruf}) was used to generate Figures 3 and 4; standard
``double precision'' floating-point operations were sufficient.  We
found, over the domain of Figure 3a, the largest difference from the
exact minimum number of pools, from Equation (\ref{rhof}),
was about six percent of the number of pools.
Also, the approximate result
typically differs from the exact result by no
more than one pool.
 
\subsection{Methods for Ranking Clones}
\label{ee}
 
Although we focus on ranking the clones according to the probability
that each is positive---given the pool-assay data and a model for the
experimental errors---ranking sets of candidate positive clones can
also be a desirable objective.  In the absence of errors, it is
usually possible to identify many resolved positive clones, but it can
also be important to rank the remaining clones according to the
probability that they are positive.  These objectives can be achieved
in the following framework.
 
Bayes' rule can be used to estimate the probability that an individual
clone is positive, given the vector of pool-assay outcomes, denoted
$V$, and a model for the experimental errors.  In this context, Bayes'
rule can be written as follows:
 
\begin{equation}
\P(I_i^+|V)=\left(1+\frac{\P(V\& I_i^-)}{\P(V\& I_i^+)}\right)^{-1},
\label{bayes}
\end{equation}
in which $I_i^+$ (respectively $I_i^-$) denotes the event that clone
$i$ is positive (respectively negative).  To evaluate exactly the
ratio on the RHS of (\ref{bayes}) would require a calculation for each
possible subset of the clones, taken to be $P$, the set of positive
clones.  Because $2^n$ is typically very large, this is not feasible,
and, therefore, we sampled subsets to estimate the ratio.  These
estimates were then used to rank the clones.  In our preliminary
studies of this approach, we selected subsets in which each clone has
an equal probability of appearing.  For the example described in
Section \ref{res} with 33,000 clones, we sampled approximately
330,000 subsets for each clone, both for the numerator and
denominator of the ratio.  To reduce both the noise of and the
computational work required for the sampling, individual clones were
added to and removed from these two groups of subsets---one used for
the numerator and the other for the denominator.  However, a much more
efficient approach is Gibbs sampling of $P$, and the Hastings-Metropolis 
algorithm have also been implemented, facilitating the estimation of Equation
(\ref{bayes}) \cite{ber}, \cite{bru}.
 
The model for experimental errors enters into the evaluation of the
RHS of Equation (\ref{bayes}).  For each set $P$ of selected clones,
taken to be the positive clones, we find the union of the pools in
which any of these clones occur.  Call this set $\U$.  In the
pool-assay data, $V$, let $v_{-|+}$ be the number of negative pools in
$\U$, and let $v_{+|-}$ be the number of positive pools not in $\U$.
Furthermore, a two parameter error model is adopted, with the
error-rate parameters for false-positive and false-negative pool
assays equal $\l_{+|-}$ and $\l_{-|+}$, respectively.  These errors
are taken to be independent; thus, the probabilities on the right-hand
side of Equation (\ref{bayes}) are evaluated using
 
    $$\P(V|P) = \l_{+|-}^{v_{+|-}}(1-\l_{+|-})^{v-|\U|-v_{+|-}}
                \l_{-|+}^{v_{-|+}}(1-\l_{-|+})^{  |\U|-v_{-|+}},$$
 
in which
$|\U|$ denotes the cardinality of $\U$.
 
\section{Results}
\label{res}
 
\subsection{Random $k$-sets}
\label{ksres}
 
In this Section we employ our criterion---the expected number of
resolved positive clones---to select optimum parameters for random
$k$-sets designs.  For libraries of varying number of clones, $n$, and
coverage, $c$, the minimum number of pools required for random
$k$-sets designs to achieve an expected number of resolved positive
clones equal to 0.5 $c$ and 0.95 $c$ is illustrated in Figures 3a and
4a, respectively.  Figures 3b and 4b depict the corresponding optimum
values of $k$---maximizing the expected number of resolved positive
clones for the number of pools plotted in Figure 3a and 4a.  The
latter plots depict a number of values of $c$ at which the optimum
value of $k$ changes abruptly.  At these values of $c$, the product of
$c$ with the probability that there are $j$ or fewer positives is
comparable to the desired expected number of resolved positive clones.
As $c$ increases through each of these ``transitions'', the optimal
pooling design resolves the cases with one additional positive clone,
in order to achieve the desired expected number of
resolved positive clones.  These transitions are most pronounced for
small values of $c$---for example, where the optimal designs go from
resolving only one positive to resolving two positives, over a small
range of $c$ centered on $0.69$.  These transitions influence the
minimum number of pools, contributing to the irregularity of the
contours in Figure 3a and Figure 4a in this region.  For most values
of the variables $c$ and $n$, the dependence of the expected number of
resolved positive clones upon $k$ is not pronounced in the vicinity of
the optimum value of $k$.  Thus, the optimum $k$ provides a rough
guideline for efficient pooling designs, as can be seen in the
following examples.
 
Given a goal of five resolved positive clones, on the average, a
tenfold-coverage, human-genome library of 33,000 clones \cite{my}
could be accommodated on 170 pools (Figure 3a). The optimum value of
$k$ is ten, which would result in an average of 1,941 clones per pool.
The expected number of unresolved negative clones is approximately 44.
If, instead, the goal were to have an an average of 9.5 resolved
positive clones, then 253 pools would be required (Figure 4a). The
optimum value of $k$ would also be ten and there would be an average
of approximately 1,304 clones per pool.  The expected number of
unresolved negative clones is approximately 2.8.  Therefore, on the
average, 3.3 confirmatory tests would be required to resolve the
status of every clone.
 
It may also be desirable to bound the average number of clones in a
pool in order to avoid high pool-assay error rates. This can be
achieved by constraining $k$, the number of pools containing any one
clone.  Suppose it were desirable to have fewer clones per pool, say,
approximately 1,000, and also to achieve an average of five resolved
positive clones for the library described in the previous paragraph.
Then, from Equation (\ref{rhof}), 191 pools would be required with
$k=6$, and the average number of clones in a pool would be
approximately 1037.
 
Experimental error necessitates more pools to achieve comparable
results.  We performed some computer-simulation experiments for a
tenfold-coverage library of 33,000 clones with a false-negative error
rate of 0.1 and a false-positive error rate of 0.01---rates consistent
with our preliminary experiments on pools containing approximately 110
clones, described in Section \ref{eufour}.  We used a random ten-sets
design, and an arbitrary set of ten clones was selected to be the
positive clones.  The result of the Bayes' ranking, described in
Section \ref{ee}, was that five of these ten clones were ranked in the
top ten.  Thus, it is feasible to identify the positive clones---even
with appreciable experimental error.
 
\subsection{Comparison with other designs}
\label{comp}
 
To facilitate comparison, we consider a tenfold-coverage, human-genome
library of 72,000 clones for which the following row-and-column
pooling design has been implemented \cite{chu}.  The library is
partitioned into 94 lots, all but one containing eight 96-well dishes
and the remaining one with six.  The rows and columns of eight
microtiter dishes are combined to construct 20 pools.  In addition,
eight more pools, each containing all of the clones from one dish, are
constructed.  Thus, the total number of pools is $28\times94=2,632$.
For this design, a lot contains a resolved positive clone only if it
contains only one positive clone.  Therefore the expected number of
resolved positive clones is approximately $10e^{-10/94} \doteq 9.0$,
and the expected number of unresolved negative clones is approximately
$3.3$ \cite{bar}.  The following designs have been proposed for
screening the same library, using approximately one-tenth as many
pools.
 
\cite{bar} proposed novel pooling designs for a
library with the aforementioned parameters.  One of these designs
assigned each clone to a lattice point in a cubic, integer lattice
with 43 points on each axis.  Each pool would contain all of the
clones with each coordinate: thus $3\times43=129$ pools would result.
To gain more information about the positive clones, a linear
transformation was used to obtain a new configuration of the clones on
the lattice points and, thus, another set of 129 pools, specified by
the new coordinates.  Thus, each cubic configuration yields three
groups of 43 pools with the property that each clone occurs in one
pool from each group.  In general, this property of the Barillot {\em
et al.} designs distinguishes them from the designs we propose.
 
Computer simulation was used for the 258-pool cubic design to
determine that the expected number of resolved positive clones is
nearly $8.8$.  Also, the expected number of unresolved negative clones
is approximately equal $13.3$. (We assumed that 72,000 of the lattice
sites were initially selected uniformly at random for the first
configuration).
 
For the same library, using the results of Section \ref{mm}, the
optimum choice of $k$ for a random $k$-sets design on 258 pools is 11.
In this case the expected number of resolved positive clones is $9.1$
and the expected number of unresolved negative clones is $5.3$.
However, if we choose $k=6$ to achieve pools comparable in size with
the ``cubic'' design, then the expected number of resolved positive
clones is $7.3$, and the expected number of unresolved negative clones
is $12.6$.  The $k$-sets packing designs, described in Section
\ref{gpd}, would yield larger expected numbers of resolved positive
clones.
 
\subsection{Screening results}
\label{eufour}

The following theoretical and computational results bear on the
predicted performance of our four-sets packing design.  As above, this
library is assumed to have a binomial number of positive clones with
expectation 2.5.  Computer simulation was used to estimate the
expected number of resolved positive clones for the four-sets packing
design at approximately 1.47---versus 1.36, from (\ref{rhof}), for
random four-sets designs.  Similarly, the expected number of
unresolved negative clones for the four-sets packing design is
3.98---versus 4.68, from (\ref{nuasymp}), for a random four-sets
design.  Thus, to identify all the positives would require
confirmatory testing of 6.48 clones, on the average.
  
The four-sets packing design was implemented for the 1298-clone, human
chromosome 16, YAC library.  We observed numerous false negative and
also false positive pool assays, precluding the identification of
clones on the basis of being either resolved positive or unresolved
positive or negative.  Twenty-two STSs were screened against the pools
to achieve closure of the chromosome 16 framework map.  We
ranked the clones according to the probability of being positive,
based on the pooling results, as described in Section \ref{ee}.
We set the error rates
$\lambda_{+|-}$ and $\lambda_{-|+}$ equal 0.14 and 0.06, respectively.
These rates are primarily based on comparing the frequency of positive
pools with that predicted using the coverage.  After performing
confirmatory testing on the top eleven clones (on average), an average
of 1.8 positive clones were identified.  Six of the twenty-two STSs
yielded no positive clones.

\section{Discussion}
\label{dis}
 
Improved methods for screening clone libraries will allow more
efficient use of currently available biological resources.  We propose
using $k$-sets designs for unique-sequence screening of large clone
libraries.  These designs are efficient, flexible, easy to specify and
can allow screening in one-pass, thus minimizing human intervention.
When possible, we also advocate using the $k$-sets packing designs,
discussed in Section \ref{gpd}, which can yield further substantial
gains in efficiency.  The automated implementation of a four-sets
packing design for a YAC library containing 1,298 clones, over a
period of 12 hours, demonstrates the utility of commercially-available
robots.
 
In addition to the expected number of unresolved negative clones,
proposed by Barillot {\em et al.,} we propose a new design performance
criterion: the expected number of resolved positive clones.
Optimizing a pooling design according to either criterion will be
sensitive to the upper tail of the distribution for the number of
positive clones, which we have assumed to be binomial for design
comparison.  For some applications, the determination of average
behavior might not be adequate; thus, it may be useful, for example,
to estimate the probability that $j$ positives are resolved positive
for a pooling design, given that there are $p$ positives.  In any
case, pooling designs which have a high probability of achieving the
screening objective are useful, even if there is no guaranteed
performance for a particular STS.  On the other hand, to guarantee
that the status of all of the clones is resolved would require many
more pools than given in Figure 4a.
 
We computed the smallest number of pools required for a random
$k$-sets design to achieve a given expected number of resolved
positive clones.  Figures 3a and 4a depict the number of pools
required for libraries of coverage $c$ with $n$ clones, and Figures 3b
and 4b depict the optimum values of $k$.  When implementing any
pooling design, it could be informative to compare the total numbers
of pools and the number of pools containing each clone to those
plotted in Figures 3 and 4 because, after optimization, these designs
come close to achieving optimal performance.  One could use the Bayes
ranking in several ways to find the best parameters for a random
$k$-sets design in the presence of experimental error.  For the time
being, we selected a design from Figure 3 and simulated its
performance in the presence of a realistic level of experimental
error.
 
Data on pool-assay errors is clearly important in the design of
pooling experiments.  Pools with a smaller proportion and
concentration of target DNA might be less likely to yield the correct
PCR product.  Also, some primers could be more prone to fail to
produce products than others.  False-positive results could result
from cross-contamination.  It is not clear that our simple model of
experimental error---with independent probabilities for false-negative
and for false-positive pool assays---is adequate.  We will assess the
fidelity of our pooling and screening experiments while assaying the
four-sets packing pools for our physical map of human chromosome 16.
 
Experimental error blurs the distinctions between clones in our four
categories and motivates both the consideration of error-correcting
pooling designs \cite{b&t} and effective ranking of the
candidate-positive clones \cite{bru}.  In general, an automated
ranking algorithm will also propose a set of clones for confirmatory
testing.  It could propose several candidate clones and use the
results from these confirmatory tests when proposing further candidate
clones.  One criterion for evaluating ranking algorithms could be the
average number of proposed clones required to identify $i$ of the
positive clones or all of them, in the event that there are fewer than
$i$ positives.  Or it could be the sum of the probabilities of being
positive assigned to positive clones in computer simulations.  One can
optimize the design of the pooling experiments, given the selection of
a ranking technique. In addition, the relative costs of confirmatory
testing and of missing positive clones could be employed as part of
design optimization.
 
It may often be possible to generate designs with better performance
than random $k$-sets designs.  Such designs could include $k$-sets
packings, in which there is a bound, $t$, on the number of pools in
which two clones coincide.  By varying the value of $t$, $k$-sets
packing designs can provide an extremely powerful and flexible
approach to library screening.  When $t=k$, we have random $k$-sets
designs which are efficient and very easy to construct.  As $t$ is
decreased, we gain even greater efficiency at the cost of additional
computations in design generation.  At the other extreme we have
maximum-size $k$-sets packings which are, we believe, maximally
efficient but often difficult to construct.  Further, we have shown
\cite{b&t} that packings which achieve the maximum possible size are
best possible in some cases, in terms of maximizing the probability of
a one-pass solution.  Such designs can also be optimal subject to
guaranteed error-detection requirements.  The magnitude of the
improvement one can achieve by constraining the intersections is
exemplified by the predicted performance of the cubic row-column
pooling design \cite{bar}.  Thus, we are exploring combinatoric
optimization techniques for the construction of $k$-sets packing
designs.  Some preliminary methods were used for optimizing the
pooling design for the Cla I, YAC library described in Section
\ref{gpd}.  We plan to improve these by applying techniques such as
the method of conditional expectations to de-randomize generation of
the designs \cite{pm}.
 
In summary, although our preliminary results should prove useful, much
exploratory work remains---based upon a better understanding of the
prevalent experimental errors---to achieve superior pooling designs,
further reducing the labor and increasing the efficiency of
large-scale, library-screening experiments.  Furthermore, pooling the
clones from a pre-existing map will involve new challenges because of
the prior knowledge about the joint probability distribution for
positive clones.
 
\begin{center}
{\bf Acknowledgements}
\end{center}
 
\noindent
This manuscript is dedicated to Drs H.\ Levine and H.M.\ McConnell of
Stanford University.  We thank Drs P.\ Hagan of Los Alamos National
Laboratory (LANL) and J.\ Percus of New York University for sharing
their insights on the asymptotic behavior of $K^{(p)}_{j}$, a
probability described in the Appendix.  We also thank Drs R. Dougherty
of Ohio State University and M.\ Goldberg of Rensselaer Polytechnic
Institute for encouraging us to consider random designs.  We thank Dr
P. Medvick of LANL for commenting on our descriptions of software.  We
thank Ms A.\ Ford of LANL for technical assistance with our first
pooling experiments, and Ms C.\ Campbell, Ms L.\ Goodwin, Ms J.\ Tesmer, and
Ms L.\ Meincke for assistance with the pooling experiments described herein.
DCT is grateful to the UK SERC for a
visiting-fellow research grant---which led to the completion of this
manuscript---and he thanks the faculty and staff of the School of
Mathematical Sciences, Queen Mary \& Westfield College, University of
London, for their hospitality.  WJB acknowledges a US Department of
Energy Distinguished Human Genome Postdoctoral Fellowship.  This work
was performed under the auspices of the US Department of Energy, and
was funded both through the Center for Human Genome Studies at Los
Alamos and by a LANL Laboratory-Directed Research and Development
grant.

\begin{center}{\bf Figure Legends}\end{center}
\noindent
Figure 1a
 
\noindent
Title:  ``Incidence Matrix for a Pooling Design''
 
A binary incidence matrix for a particular pooling scheme with six
pools and five clones.  A clone occurs in a pool if the corresponding
matrix element equals unity and does not occur in a pool if the corresponding
matrix element equals zero.  Pool number 1 contains clones number 1, 3, and 4,
pool number 4 contains only clone number 5, et cetera.\medskip
 
\noindent
Figure 1b
 
\noindent
Title:  ``Terminology and Example of Categorization''
 
The binary incidence matrix of Figure 1a is depicted, but now clones 1
and 2 are taken to be positive clones and clones 3, 4 and 5 are taken
to be negative.  The `1's for the positive clones are replaced by
`+'s, the `1's for the negative clones are replaced by `-'s, and the
`0's are omitted.  In the absence of experimental error, all the pools
containing either of the positive clones would be positive; thus, only pool 4
would be negative, as depicted above the horizontal line opposite `Pool Assay'.
This hypothetical pooling experiment would resolve the status of clones 1 and
5, and the remaining clones could be either positive or negative.
Clone 5 is a negative clone occurring in pool 4,
which contains no positive clone, thus it is a resolved negative
clone.  Clones 3 and 4 are negative clones occurring only in pools containing positive
clones; thus, they are both unresolved negative clones.  Because clone 1 is a positive
clone occurring only in
pool 6,  which does not contain either another positive clone or an unresolved
negative clone, it is a resolved positive clone.  Clone 2 is a positive clone occurring only
in pools containing either other positive clones or unresolved negative clones; thus,
it is unresolved positive.\medskip
 
\noindent
Figure 2
 
\noindent
Title:  ``Row-and-Column Pools of the Clones in a 96-well dish''
 
A row and column design with 96 clones and 20 pools. Two clones are
positive (row 3, column 3 and row 5, column 8) and hence the four
indicated pools would be positive, in the absence of experimental error.
There are two unresolved negative clones (row 5, column 3 and row 3, column 8)
and 92 resolved negative clones.
The two positive clones are unresolved positive.\medskip
 
\noindent
Figure 3a
 
\noindent
Title for Figure 3a: ``Minimum Number of Pools with the Expected
Number of Resolved Positives Equal $0.5c$''
 
The $x$ axis is $n$, the number of clones to be
pooled; $1000 \le n \le 100,000$, and the $y$ axis is $c$,
the coverage parameter; $ 1/4 \le c \le 16$.  Both axes have a logarithmic scale,
and the tics facing the inside of the plot are at integer values while those on
the outside are at values corresponding to half integers.
The smallest value of $v$ such that a random $k$-sets
design can achieve the target expected number of 0.5$c$
resolved positive clones is depicted.  Equation (\ref{rhoruf})
was used to generate the data for this plot.\smallskip
 
\noindent
Figure 3b
 
\noindent
Title for Figure 3b: ``Optimum $k$, Expected Number of Resolved
Positives Equal $0.5c$''

The $x$ axis is $n$, the number of clones to be
pooled; $1000 \le n \le 100,000$, and the $y$ axis is $c$,
the coverage parameter; $ 1/4 \le c \le 16$.  Both axes have a logarithmic scale,
and the tics facing the inside of the plot are at integer values while those on
the outside are at values corresponding to half-integers.
The values of $k$ achieving the maximum expected
number of resolved positive clones are depicted, with the value of $v$
depicted in Figure 3a. These expected values slightly exceed $0.5 c$.
For this domain of $c$ and $n$, optimal values of $k$ fall between 6
and 12.  Equation (\ref{rhoruf})
was used to generate the data for this plot.\medskip
 
\noindent
Figure 4a
 
\noindent
Title for Figure 4a: ``Minimum Number of Pools with the Expected
Number of Resolved Positives Equal $0.95c$''

The $x$ axis is $n$, the number of clones to be
pooled; $1000 \le n \le 100,000$, and the $y$ axis is $c$,
the coverage parameter; $ 1/4 \le c \le 16$.  Both axes have a logarithmic scale,
and the tics facing the inside of the plot are at integer values while those on
the outside are at values corresponding to half-integers.
The smallest value of $v$ such that a random $k$-sets
design can achieve the target expected number of 0.95$c$
resolved positive clones is depicted. Equation (\ref{rhoruf})
was used to generate the data for this plot.\smallskip
 
\noindent
Figure 4b
 
\noindent
Title for Figure 4b: ``Optimum $k$; Expected Number of
Resolved Positive Clones Equal to $0.95c$''
 
The $x$ axis is $n$, the number of clones to be
pooled; $1000 \le n \le 100,000$, and the $y$ axis is $c$,
the coverage parameter; $ 1/4 \le c \le 16$.  Both axes have a logarithmic scale,
and the tics facing the inside of the plot are at integer values while those on
the outside are at values corresponding to half-integers.
The values of $k$ achieving the maximum expected
number of resolved positive clones are depicted, with the value of $v$
depicted in Figure 4a. These expected values slightly exceed $0.95
c$.  For this domain of $c$ and $n$, optimal values of $k$ fall
between 7 and 12.  Equation (\ref{rhoruf})
was used to generate the data for this plot.\bigskip
 
\appendix
\section{Appendix:  Derivations}
\label{app}
 
The inner summation in (\ref{nuexact}), which equals the probability,
$K^{(p)}_k$, that $k$ specified pools contain one or more
positive clones when there are
exactly $p$ positive clones, can be evaluated via the recursive
formula
\begin{equation}
 K^{(p)}_{j}=\sum_{i=0}^{j}{j\ch i}{v{-}j\ch k{-}i}
{v\ch k}^{-1}K^{(p-1)}_{j-i},
\label{covjr}
\end{equation}
for $j\le k$, where $K^{(0)}_{j-i}=1$ if $j=i$, otherwise
$K^{(0)}_{j-i}=0$.  Equation (\ref{covjr}) is numerically advantageous
because it involves no subtractions.
 
We turn now to the expected number of unresolved positive clones,
$\bar P$. Let $\a^{(p)}$ denote the probability that a selected
positive clone is unresolved, given that there are exactly $p$
positive clones.  When $p=1$, the positive clone will be unresolved
only if a negative clone occurs in precisely the same $k$ pools, so
that
\begin{equation}
\a^{(1)} = 1-\left(1-{v\ch k}^{-1}\right)^{n-1}.
\label{p1}
\end{equation}
(In practice, a $k$-sets design would usually be
generated so that no two clones occupy precisely the same pools and hence
$\a^{(1)}=0$, as described in Section \ref{gpd}.  Here, however,
it is convenient to consider standard random $k$-sets designs.)
 
For $p\ge 2$, we determine $\a^{(p)}$ by conditioning on the values of
three random variables $U$, $X$, and $Y$, where $U$ is the number of
unresolved negative clones, $X$ is the number of pools which would be
positive if the selected positive clone were removed, and $Y$ is the
number of positive pools.  Thus, $Y{-}X$ is the number of pools
containing the selected positive clone but no other positive clone.
The selected positive clone will be unresolved positive either if
$Y{-}X$ is zero or if each of these $Y{-}X$ pools contains at least
one unresolved negative clone.  Thus,
\begin{equation}
\a^{(p)}=\sum_{u=0}^{n-p}\sum_{x=k}^v\sum_{y=x}^v\P(A|U{=}u,X{=}x,Y{=}y)
\P(U{=}u,X{=}x,Y{=}y),
\label{prp}
\end{equation}
where $A$ denotes the event that every pool in which the selected
positive clone occurs contains either a positive clone or an
unresolved negative clone.  There is an implicit conditioning on the
number $p$ of positive clones in each term of (\ref{prp}).
 
To evaluate (\ref{prp}), we use the equality
$$
\P(U{=}u,X{=}x,Y{=}y)=\P(U{=}u|X{=}x,Y{=}y)\P(Y{=}y|X{=}x)\P(X{=}x).$$
Now $$\P(X{=}x)={v\ch x}L^{(p-1)}_{v-x},$$ where $L^{(p)}_j$
denotes the probability that $j$ specified pools are precisely the
negative pools.  By the inclusion-exclusion principle,
\begin{equation}
L^{(p)}_j=\sum^{v-k}_{i=j}{v{-}j\ch i{-}j}(-1)^{i-j} z_{i}^{p},
\label{uncovj}
\end{equation}
in which $z_i$ is the probability, introduced at (\ref{nuexact}), that
a given clone occurs in none of $i$ specified pools. Therefore
\begin{equation}
\P(X{=}x)={v\ch x}\sum_{i=v-x}^{v-k}{x\ch v{-}i}(-1)^{i-v+x}z_i^{p-1}.
\label{q}
\end{equation}
As was the case for $K^{(p)}_j$, a recursive formula for $L^{(p)}_j$
is also available:
\begin{equation}
L^{(p)}_j=\sum_{i=j}^{j+k}{i\ch j} {v{-}i\ch
k{-}i{+}j}{v\ch k}^{-1}L^{(p-1)}_{i},
\label{uncovjr}
\end{equation}
in which $L^{(0)}_{i}=1$ if $i=v$, otherwise $L^{(0)}_{i}=0$.
Parenthetically, the inclusion-exclusion principle can be used
to derive the $K_{j}^{(p)}$ from the $L_{j}^{(p)}$ and vice-versa.
 
          $$K_{j}^{(p)} =\sum_{i=0}^{v-j}{v-i\ch i} L_{i}^{(p)};$$
 
$$L_{j}^{(p)} =\sum_{i=v-j}^{v}{j\ch i-(v-j)}(-1)^{i-(v-j)} K_{i}^{(p)}.$$
 
Given $X{=}x$, the distribution of $Y$ is Hypergeometric:
\begin{equation}
\P(Y{=}y|X{=}x)={x\ch y{-}k}{v{-}x \ch y{-}x}{v
\ch k}^{-1}.
\label{r}
\end{equation}
Further, given $Y{=}y$, the distribution of $U$ is conditionally
independent of $X$ and is binomial:
\begin{equation}
\P(U{=}u|Y{=}y) = B(n{-}p,u,\b),
\label{s}
\end{equation}
where $B(a,b,t)$ denotes the binomial probabilities (\ref{bin}) and
where $\b$ is the probability that a given negative clone is
unresolved negative, so that $$\b={y\ch k}\,{v\ch k}^{-1}.$$
Finally, for $u \ge 0$,
\begin{equation}
\P(A|U{=}u,X{=}x,Y{=}y)=\sum_{j=0}^{y-x}{y{-}x\ch j}(-1)^j\xi_j^u,
\label{prguxyp}
\end{equation}
in which $$\xi_j = {y{-}j\ch k}\,{y\ch k}^{-1}.$$ The RHS of
(\ref{prguxyp}) is, essentially, $K_{y-x}^{(u)}$, noting that each of
the possible ${y\ch k}$ subsets of the $y$ positive pools has equal
probability of being the pools which contain a particular unresolved negative clone.
 
Equation (\ref{rhof}), follows from (\ref{p1}) and (\ref{prp}),
together with (\ref{q}), (\ref{r}), (\ref{s}), and (\ref{prguxyp});
\begin{equation}
\bar{P}=\sum_{p=1}^npB(n,p,c/n)\a^{(p)},
\label{last}
\end{equation}
where $B(a,b,t)$ denotes the binomial probabilities (\ref{bin}).

\end{document}